\documentclass[graybox]{svmult}

\usepackage{type1cm}        % activate if the above 3 fonts are
\usepackage{makeidx}         % allows index generation
\usepackage{graphicx}        % standard LaTeX graphics tool
\usepackage{multicol}        % used for the two-column index
\usepackage[bottom]{footmisc}% places footnotes at page bottom

\usepackage{newtxtext}       % 
\usepackage[varvw]{newtxmath}       % selects Times Roman as basic font

\usepackage{amsmath,amsthm}
\usepackage{amsfonts}
\usepackage{mathtools}
\usepackage{bbm}
\usepackage{mathrsfs}

\usepackage{tikz}\usetikzlibrary{decorations.markings,decorations.pathreplacing, patterns, calc,math}
\usepackage{circuitikz}
\usepackage{cite}
\usepackage{hyperref}

\newcommand{\R}{\mathbb{R}}
\newcommand{\N}{\mathbb{N}}
\newcommand{\e}{\varepsilon}
\newcommand{\pa}{\partial}
\newcommand{\dl}{\delta}
\newcommand{\fb}{\beta}
\newcommand{\fbv}{\beta^{\mathrm{var}}}
\newcommand{\BC}{\mathcal{M}}
\newcommand{\SC}{\mathcal{S}}
\newcommand{\uBC}{\overline{\mathcal{M}}}
\newcommand{\lBC}{\underline{\mathcal{M}}}
\newcommand{\uSC}{\overline{\mathcal{S}}}

\newcommand{\ubex}{\overline{\mathfrak{m}}}
\newcommand{\lbex}{\underline{\mathfrak{m}}}
\newcommand{\bex}{\mathfrak{m}}
\newcommand{\usex}{\overline{\mathfrak{s}}}
\newcommand{\lsex}{\underline{\mathfrak{s}}}

\newcommand{\udimout}{\overline{\dim}^{\mathrm{out}}}
\newcommand{\ldimout}{\underline{\dim}^{\mathrm{out}}}
\newcommand{\dimout}{\dim^{\mathrm{out}}}
\newcommand{\Minkout}[1]{\mathcal{M}^{{\rm out},#1}}
\newcommand{\uMinkout}[1]{\overline{\mathcal{M}}^{{\rm out},#1}}
\newcommand{\lMinkout}[1]{\underline{\mathcal{M}}^{{\rm out},#1}}
\newcommand{\sH}{\mathscr{H}}
\newcommand{\1}{\mathbbm{1}}
\newcommand{\dist}{\mathrm{dist}}

\newcommand{\Unp}{\mathrm{Unp}}

\renewcommand{\o}{\omega}
\newcommand{\z}{\zeta}

\renewcommand{\E}{\mathrm{e}}

\renewcommand{\I}{\mathbbm{i}}

\renewcommand{\d}{\mathrm{d}}

\newcommand{\udim}{\ov{\dim}}
\newcommand{\ldim}{\underline{\dim}}
\newcommand{\uS}{\overline{\mathcal{S}}}

\newcommand{\C}{\mathbb{C}}
\newcommand{\eS}{\mathbb{S}}

\newcommand{\eps}{\varepsilon}	 % Epsilon

\newcommand{\frk}{\mathfrak}

\newcommand{\ov}[1]{\overline{#1}}
\newcommand{\un}[1]{\underline{#1}}

\newcommand{\re}{\operatorname{Re}}

\newcommand{\Z}{\mathbb{Z}}

\newcommand{\res}{\operatorname{res}}

\def\sideremark#1{\ifvmode\leavevmode\fi\vadjust{\vbox to0pt{\vss % the remark
    \hbox to 0pt{\hskip\hsize\hskip1em           %                will appear only
 \vbox{\hsize3.5cm\tiny\raggedright\pretolerance10000%                on the side
 \noindent #1\hfill}\hss}\vbox to8pt{\vfil}\vss}}}%
\newcommand\equidistantPoints[4]{%
    \draw[#1] (0, 0) rectangle (#2, #2);
    \pgfmathsetmacro{\scale}{#4}
    \pgfmathsetmacro{\step}{#2 / (#3 + 1)}
    \pgfmathsetmacro{\last}{#3-1}
    \foreach \i in {0,...,\last} {
      \foreach \j in {0,...,\last} {
        \fill[#1] (\i*\step+\step, \j*\step+\step) circle [radius=0.03/\scale];
      }
    }
}

\makeindex             % used for the subject index
\begin{document}

\title*{Review of Steiner formulas in Fractal Geometry via Support measures and Complex Dimensions}
\titlerunning{A Review of Steiner formulas in Fractal Geometry} 
\author{Goran Radunovi\'c\orcidID{0000-0002-4791-9212}}
\institute{Goran Radunovi\'c \at University of Zagreb, Faculty of Science, Department of Mathematics, Horvatovac 102a, 10000 Zagreb, Croatia, \email{goran.radunovic@math.hr}}
\maketitle

\abstract{We review the theoretical framework that establishes a crucial bridge between the general Steiner-type formula of Hug, Last, and Weil and the theory of complex (fractal) dimensions of Lapidus et all. Two novel families of geometric functionals are introduced based on the support measures of the set itself as well as of its parallel sets, respectively. The associated scaling exponents provide new tools for extracting geometric information of the set, beyond its classical fractal
dimensions while also encoding its outer Minkowski dimension.  Furthermore the scaling exponents also directly connect to the complex  dimensions of the set while preserving essential geometric information. The framework provides a fundamental link between measure-theoretic approaches and analytical methods in fractal geometry, offering new perspectives on both the geometric measure theory of singular sets and the complex analytic theory of fractal zeta functions.}

\section{Introduction and Historical Context}

The study of fractal geometry has undergone significant developments over the past several decades, with researchers continuously seeking to extend classical geometric concepts to increasingly singular and irregular sets.
While fractal dimension (Hausdorff, box-counting, Minkowski, or other) provides valuable insights into scaling behavior and space-filling capacity, it cannot fully characterize the complex geometric structure of highly irregular sets and it does not encode enough information in order to be be useful in some applications.
Many fractal structures can share identical scaling relationships while exhibiting dramatically different geometric appearance, the Devil's staircase (graph of the Cantor function) being probably the most famous example having all common fractal dimensions as well as its topological dimension equal to 1. This limitation has motivated researchers to develop additional geometric functionals that provide more nuanced information about fractal sets.

Already Mandelbrot has proposed the notion of `lacunarity' (directly connected to the Minkowski content) \cite{Ma94-lacun} to describe finer detail of fractal sets.
It has found many applications in other fields of mathematics such as dynamical systems, e.g.,  \cite{Res2013,Res2014,ZupZub}, spectral theory \cite{LaPom,CaBr86,LaPo96}, number theory \cite{LapHer21,LaPom},  applied mathematics, e.g., \cite{Ba13,DEMES13,DiIeva24}, etc.

Beyond Minkowski content, several geometric descriptors have been proposed and applied in the applied mathematics literature to extract additional geometric information from fractal structures, see e.g.\ \cite{bruno2012shape,CaGa13,DEMES13}. However, many of these concepts lack rigorous mathematical foundations. The development and use of these alternative approaches demonstrates the critical need for more comprehensive geometric descriptors that can adequately characterize the complex properties of fractal sets.

Several more sophisticated mathematical tools have already emerged to address this limitation.
Building upon Federer's classical curvature measures \cite{Fed59}, also known as intrinsic volumes or  Minkowski functionals, researchers have developed {\em fractal curvature measures} extending the notion of curvature to highly irregular sets~\cite{winter2014minkowski, Wi08, Zahle2011,spodarev2014estimation, klatt2020geometric}. These measures arise as weak limits of appropriately rescaled curvature measures of parallel sets and provide geometric information beyond just the fractal  dimension~\cite{lapidus2018fractal}.

The {\em theory of complex dimensions}, developed through the use of fractal zeta functions, provides another powerful framework for characterizing fractality~\cite{FZF,LF12,lapidus2017distance,lapidus2018fractal}. These complex-valued functions encode scaling properties and oscillatory behavior that single real-valued dimensions cannot capture, offering refined geometric information through their pole structure.

The seminal work of Hug, Last and Weil \cite{HugLasWeil}, established a comprehensive general Steiner-type formula for arbitrary closed sets by introducing the support measures. This represented a significant non-additive generalization of Federer's classical curvature measures \cite{Federer,Fed59}, extending their applicability far beyond sets with positive reach to general closed sets in Euclidean space.
Their work emerged from a rich historical context that traces back to the classical works of Steiner and his fundamental formula for convex bodies. Based on earlier work of Stach\'o \cite{Stacho79}, Steiner's original formula \cite{Sta} expressed the volume of parallel sets of convex bodies, i.e., their $\eps$-neighborhoods,  as polynomials in the dilation parameter $\eps$, with coefficients given by the intrinsic volumes. This elegant relationship between geometric dilation and algebraic structure provided a foundation for understanding how geometric objects scale under parallel expansion.

However, while the Hug-Last-Weil theory provides a mathematically rigorous and general framework and has found strong application in stochastic geometry \cite{HugLast00,HugLasWeil,Vi25}, and geometry measure theory \cite{HugSan2022}, there are significant challenges to applying it to fractal sets. The support measures, though well-defined in complete generality, do not clearly reveal the classical fractal properties that have become central to the field, such as Minkowski dimension and content, or the complex fractal dimensions \cite{FZF,LF12}.

This gap between the measure-theoretic generality and the specific needs of fractal analysis has motivated the author and Steffen Winter in \cite{RaWi2,RaWi1}  to offer a significant step toward bridging this gap. The contribution is particularly important because it demonstrates how the general Steiner formula can be systematically decomposed into components that separately capture different geometric aspects of fractal sets, enabling a more refined and comprehensive analysis of fractal properties than was previously possible while simultaneously connecting the Hug-Last-Weil's theory directly to the theory of complex dimensions of M.\ L.\ Lapidus et all.

Note also that recently support measures have been generalized for sets in non-Euclidean spaces in \cite{HugSan2022}. 
Hence, the scaling exponents and contents reviewed here can also be naturally generalized to such settings.
However, we restrict our review to the Euclidean setting.

% The first direction involves the theory of fractal curvature measures, which seeks to extend classical curvature concepts to fractal sets while maintaining additive properties that are natural for geometric measure theory. This approach has proven particularly valuable for understanding how fractal sets can be decomposed into simpler components, with the curvature measures providing information about the geometric complexity of these decompositions.

% The second direction encompasses the theory of complex dimensions via fractal zeta functions, pioneered by Lapidus and van Frankenhuijsen \cite{LF12} and continued by numerous collaborators \cite{FZF}.... This approach employs techniques from complex analysis to study the scaling properties of fractal sets, revealing deep connections between geometric properties and analytic structures. The complex dimensions, defined as poles of meromorphically continued fractal zeta functions, provide detailed information about the oscillatory behavior and scaling properties of fractal sets.

\section{Support measures and the general Steiner formula}

In this section we recall the framework provided by Hug, Last, and Weil \cite{HugLasWeil} based on the previous work of Stacho \cite{Stacho79}. 
Throughout this discussion, we work with nonempty compact sets $A \subseteq \R^d$ where $d \geq 2$, and we define the index set $I_d := \{0, \ldots, d-1\}$ that will appear repeatedly throughout the paper.

The geometric foundation of the theory uses the notion of parallel sets. For any compact set $A \subseteq \R^d$ and any $\e > 0$, the closed $\e$-parallel set (or $\e$-neighborhood) is defined as $A_\e := \{x \in \R^d : \dist(x,A) \leq \e\}$, representing the set of all points whose Euclidean distance to $A$ is at most $\e$. 

The key insight of the Hug, Last and Weil approach lies in the introduction of the generalized normal bundle $N(A)$ for an arbitrary closed set $A$. This construction requires careful consideration of the metric projection properties of the set. Let $\Unp(A)$ denote the set of all points $y \in \R^d \setminus A$ that possess a {\bf unique} nearest point in $A$, with $\pi_A(y)$ denoting this metric projection.
Recall also the the notion of the {\em directional metric projection}
$\Pi_A:\Unp(A)\setminus A\to N(A)$ defined as
\begin{equation}\nonumber
    \Pi_A(y):=\left(\pi_A(y),\frac{y-\pi_A(y)}{|y-\pi_A(y)|}\right)=\left(\pi_A(y),\frac{y-\pi_A(y)}{\dist(y,A)}\right).
\end{equation}
The generalized normal bundle is then defined as
$$N(A) := \left\{\Pi_A(y) : y \in \Unp(A) \setminus A\right\}\subseteq\partial A\times\mathbb{S}^{d-1}.$$

This construction extends the classical normal bundle concept from smooth manifolds to arbitrary closed sets, providing a measure-theoretic framework for understanding the ``normal directions'' associated with singular sets. Each element $(x,u) \in N(A)$ consists of a footpoint $x \in \partial A$ and a unit vector $u \in S^{d-1}$ representing an outward normal direction.

Central to the theory is also the local reach function $\dl(A, \cdot): N(A) \to [0,\infty]$, defined by $\dl(A,x,u) := \sup\{t \geq 0 : \pi_A(x+tu) = x\}$. This function measures the maximal distance one can travel from a footpoint $x\in\pa A$ in the normal direction $u$ while still maintaining $x$ as the unique nearest point in $A$. For convex sets, the local reach is always infinite, reflecting the global nature of convexity, while for general sets, it provides crucial information about the local geometric structure.
The {\em reach} of a nonempty set $A\subseteq \R^d$ is then defined as the infimmum of the local reach over all points in $N(A)$, and if it is nonzero, then the set $A$ is said to have {\em positive reach}.

The fundamental achievement of \cite{HugLasWeil} was to establish the existence of signed measures $\mu_0(A;\cdot), \ldots, \mu_{d-1}(A;\cdot)$ on the normal bundle $N(A)$, called {\em support measures}, that satisfy a general Steiner formula stated in the following theorem.\footnote{For any signed measure $\mu$, we denote by  $|\mu|$ the total variation measure of $\mu$.}

\begin{theorem}[General Steiner formula {\cite[Thm.~2.1]{HugLasWeil}}]
	\label{thm:generalSteiner}
  For any nonempty closed set $A\subseteq\R^d$, there exist signed measures $\mu_0(A;\cdot),\ldots,\mu_{d-1}(A;\cdot)$ on $N(A)$, called the \emph{support measures} of $A$, satisfying
  \begin{align}
   \int_{N(A)}\1\{x\in B\}\min\{r,\dl(A,x,u)\}^{d-i}|\mu_i|(A;\d(x,u))<\infty,
  \end{align}
  for $i=0,\ldots,d-1$, any compact set $B\subseteq\R^d$ and any $r>0$, such  that, for any measurable function $f:\R^d\to\R$ with compact support, one has that
   \begin{align}
	\label{eq:generalSteiner}
\nonumber\int_{\R^d\setminus A} & f(x) \d x \\
&=\sum_{i=0}^{d-1}\o_{d-i}\int_0^{\infty}
\int_{N(A)}t^{d-i-1}\1{\{t<\dl(A,x,u)\}}
f(x+tu)\mu_i(A;\d(x,u))\d t,
\end{align}
where $\o_k:=\frac{k\pi^{k/2}}{\Gamma(k/2+1)}$ is the surface area of the $k$-dimensional unit ball in $\R^k$.
\end{theorem}

Choosing $f=\1_{A_\e}$ in \eqref{eq:generalSteiner}, one obtains an expression for the $\e$-parallel volume of $A$ for any $\e>0$, see also \cite[\S4.3]{HugLasWeil}, namely
\begin{equation}
	\label{eq:v-par}
	V(A_\e\setminus A)=\sum_{i=0}^{d-1}\o_{d-i}\int_0^{\e}t^{d-i-1}\int_{N(A)}\1{\{t<\dl(A,x,u)\}}\mu_i(A;\d(x,u))\d t,
\end{equation}
where here, and throughout, $V(C):=\mathscr{L}^d(C)$ denotes the $d$-dimensional volume or Lebesgue measure of a set $C\subseteq\R^d$.

Theorem \ref{thm:generalSteiner} represents a far-reaching generalization of classical Steiner formula, extending its applicability from the realm of convex geometry to the full generality of closed sets in Euclidean space. The support measures $\mu_i(A;\cdot)$ encode comprehensive geometric information about the set $A$, generalizing the classical curvature measures while maintaining the essential structure needed for geometric analysis.
%The have found useful applications in stochastic geometry \cite{HugLast00,Last06} and geometric measure theory \cite{HugSan2022}.
%However, they have not been much used in fractal geometry, which we aim to remedy in \cite{RaWi1,RaWi2}.

The support measures admit an explicit integral representation that expresses them directly in terms of the generalized principal curvatures of the set \cite[Corollary 2.5]{HugLasWeil}.
Namely, for $\sH^{d-1}$-almost every point $(x,u) \in N(A)$, one can define {\em generalized principal curvatures} $k_1(A,x,u), \ldots, k_{d-1}(A,x,u)$ that extend the classical differential geometric concept to arbitrary closed sets.
Here and throughout, $\sH^{d-1}$ denotes the $(d-1)$-dimensional Hausdorff measure on the (generalized) normal bundle $N(A)$.\footnote{Appropriately normalized so that $\mathscr{H}^k([0,1]^k)=1$ for any integer $k$.}
Note that they were initially defined for sets of positive reach in \cite{zahle86} (see also \cite[\S 4.4]{RZ19}) and then the definition was extended in \cite{HugLasWeil} to the current generality. These principal curvatures determine the elementary symmetric functions $H_j(A,x,u)$ that appear in the explicit representation
\begin{equation}
    \mu_i(A; \cdot) = \frac{1}{\o_{d-i}} \int_{N(A)\cap\,\cdot} H_{d-1-i}(A,x,u) \sH^{d-1}(d(x,u)).
\end{equation}
Furthermore, for any $j\in I_d$, the elementary symmetric function $H_{j}$ of $A$ is defined (for $\sH^{d-1}$-almost all $(x,u)\in N(A)$) by
\begin{align} \label{eq:Hj}
    H_j(A,x,u):=\frac{\sum_{|I|=j}\prod_{l\in I}k_l(A,x,u)}{\prod_{i=1}^{d-1}\sqrt{1+k_i(A,x,u)^2}},
\end{align}
see e.g.\ \cite[Eq.\ (2.13)]{HugLasWeil}.
%denoted by $k_i(A,x,u)$ for $i=1,\ldots,d-1$.
%Note that the generalized principal curvatures

We point out several fundamental characteristics of support measures, for further details see \cite[\S 4.4]{HugLasWeil}. These measures are \emph{locally defined}, i.e., whenever $A_1\cap U=A_2\cap U$ for some closed sets $A_1,A_2\subseteq\R^d$ and some open set $U\subseteq\R^d$, the relation
$
\mu_i(A_1;B)=\mu_i(A_2;B)
$
is satisfied for all Borel sets $B\subseteq U\times \eS^{d-1}$ and any $i\in I_d$.

Furthermore, support measures are \emph{motion covariant}: given any $i\in I_d$, any closed set $A\subseteq\R^d$ and any rigid motion $g\in\R^d$ with rotational part $\hat g\in SO_d$ and translational part $b\in\R^d$ (so that $g(x)=\hat g(x)+b$, $x\in\R^d$), we have
$
\mu_i(g(A); g(B))=\mu_i(A;B)
$
for any Borel set $B\subseteq \R^d\times \eS^{d-1}$, where $g(B):={(g(x),\hat g(u)): (x,u)\in B}$. 
Additionally, for any $j\in I_d$ the $j$-th support measure is \emph{homogeneous of degree $j$}: given any closed set $A\subseteq\R^d$ and any scaling factor $\lambda>0$, the identity
\begin{align} \label{eq:homogen}
\mu_j(\lambda A; \lambda B)=\lambda^j \mu_j(A;B),
\end{align}
is valid for any Borel set $B\subseteq \R^d\times \eS^{d-1}$, where $\lambda B:={(\lambda x,u): (x,u)\in B}$.

\section{Fractal Curvature Measures}

In this section we review certain relationships between support measures and (generalized) curvature measures. For any compact set $A\subseteq\R^d$ with positive reach, curvature measures $C_0(A,\cdot),\ldots, C_{d-1}(A,\cdot)$ are defined via the local Steiner formula established by Federer \cite{Fed59} and they match the corresponding support measures, that is, $C_i(A,\cdot)=\mu_i(A;\cdot)$ for $i\in I_d$.
Although, support measures represent a generalization of Federer's curvature measures to arbitrary closed sets, it is important to keep in mind that they constitute a \emph{non-additive} extension of these curvature measures. 
Other (additive) extensions of the curvature measures exist.
In particular, for the class of UPR-sets, i.e., the class of sets that can be expressed as locally finite unions of sets with positive reach such that their finite intersections also possess positive reach. 
A subclass is called the {\em extended convex ring} ECR, which includes all subsets of $\R^d$ that can locally be expressed through finite unions of compact convex sets. 
Curvature measures can be extended additively to UPR-sets, implying that these measures are well-defined for any such set and satisfy additivity:
\begin{align} \label{eq:curv:additive}
C_i(A_1\cup A_2, \cdot)=C_i(A_1,\cdot)+C_i(A_2,\cdot)-C_i(A_1\cap A_2,\cdot)
\end{align}
whenever $A_1,A_2,A_1\cup A_2,A_1\cap A_2$ are UPR-sets, see \cite[\S 5.2, in particular, Cor.~5.15]{RZ19}.

Generalized curvature measures also have their support on a subset of $\partial A\times \eS^{d-1}$ but generally their support may be strictly larger than $N(A)$. However, if $A\subseteq\R^d$ has positive reach, then $\mu_i(A;\cdot)=C_i(A,\cdot)$ for any $i\in I_d$. The general relationship between these measures for the class of UPR sets is the following. For an arbitrary UPR-set $A\subseteq\R^d$ and any $i\in I_d$,
$
\mu_i(A;\cdot)=C_i(A;N(A)\cap\cdot).
$
For $A$ in the extended convex ring this is established in ~\cite[eq.~(3.1)]{HugLasWeil}, see also \cite[Thm.~3.3]{HugLast00}. In the general case, it follows by comparing the integral representations of $\mu_i(A;\cdot)$ and $C_i(A,\cdot)$ (see e.g.~\cite[Thm.~5.18]{RZ19}).
Consequently, in general, curvature measures may capture strictly more geometric information about a set $A$ than support measures by capturing how the UPR set $A\subseteq \R^d$ is represented (as a union) via an appropriate {\em index function}. On the other hand the support measures only see the set $A$ itself and the information of non-convexity of the set $A$ is directly encoded in the reach function that appears in formulas involving the support measures.

Note that generalized curvature measures (when they exist) share with the corresponding support measures the properties of being locally defined, motion covariant, and homogeneous.
Furthemrore, for any compact set $A\subseteq\R^d$ for which curvature measures exist, the total masses
$
C_i(A):=C_i(A,\R^d\times\eS^{d-1})
$
of the curvature measures $C_i(A,\cdot)$, $i\in I_d$, are called the \emph{total (Lipschitz-Killing) curvatures} of $A$ or, within convex geometry, \emph{intrinsic volumes}. For any compact set $A\subseteq\R^d$ with positive reach, we have $C_i(A)=\mu_i(A;\R^d\times\eS^{d-1})$.

\section{Basic content and associated scaling exponents}

Our approach in \cite{RaWi1,RaWi2} to making the support measures more usable in fractal geometry as well as  building a bridge between the Steiner formula of Hug, Last and Weil \cite{HugLasWeil} and the complex dimensions theory of Lapidus et all \cite{FZF,LF12}, lies in the systematic decomposition of the general Steiner formula into components that directly reveal fractal properties. This decomposition is achieved through the introduction of geometric functionals called the {\em basic contents} of closed sets arising from the associated {\em basic functions}, which serve as the fundamental building blocks for extracting dimensional information from support measures.
Later they also become a fundamental bridge connecting with the theory of fractal zeta functions and complex dimensions as we explain in Section \ref{sec:compdim} below. 

In \cite{RaWi1} we introduce for any compact set $A \subseteq \R^d$ and any index $i \in I_d$, the $i$-th basic function as
\begin{equation}
	\fb_i(t) := \fb_i(A;t) := \int_{N(A)} \1\{t < \dl(A,x,u)\} \mu_i(A; d(x,u)), \quad t > 0.
\end{equation}

Rather than working directly with parallel volumes or surface areas, the basic functions isolate the contributions of individual support measures to the overall geometric structure. This isolation proves crucial for understanding how different geometric features contribute to the fractal properties of the set.

The basic functions possess several important properties that make them particularly suitable for fractal analysis. The total variation analog 
\begin{equation}
\fbv_i(t) := \int_{N(A)} \1\{t < \dl(A,x,u)\} |\mu_i|(A; d(x,u))
\end{equation}
 is nonnegative, monotonically decreasing, and right continuous.

The key result about these basic functions is that their asymptotic behavior as $t \to 0^+$ directly encodes information about the outer Minkowski dimension of the set.
Recall that, for $q\geq 0$, the $q$-dimensional upper \emph{outer Minkowski content} of a compact set $A\subseteq\R^d$ is defined by
\begin{equation}
        \label{eq:outer-Mink}
		\uMinkout{q}(A):=\limsup_{\eps\to 0^+}\eps^{q-d} V(A_\eps\setminus A)
	\end{equation}
and the lower \emph{outer Minkowski content} $\lMinkout{q}(A)$ similarly replacing $\limsup$ by $\liminf$. Moreover, the \emph{upper outer Minkowski dimension} is then given naturally by
\begin{equation}
    \udimout_M(A):=\inf\{q\geq 0: \uMinkout{q}(A)=0\}
\end{equation}
and a lower counterpart $\ldimout_M(A)$ similarly using $\lMinkout{q}(A)$.
Although the notion of classical Minkowski content and dimension is usually studied, for our purposes the outer versions are just the right notion in general.
Note that that for compact sets $A\subseteq \R^d$ with $V(A)=0$, the outer Minkowski content and dimension coincide with the ordinary ones.

The study of asymptotics of the basic functions leads then naturally to the definition of basic contents and basic scaling exponents. Namely, the {\em $q$-dimensional upper $i$-th basic content} of a compact set $A$ is defined as
\begin{equation}
	\uBC_i^q(A) := \limsup_{t \to 0^+} t^{q-i} \fb_i(t),
\end{equation}
with corresponding total variation analog denoted by $\uBC_i^{\mathrm{var},q}(A)$ and also lower (total variation) variants (denoted by an underline). These contents measure the rate at which the basic functions scale as the parameter $t$ approaches zero, providing direct information about the dimensional properties of the set. Then, the {\em upper $i$-th basic scaling exponent} is defined as
\begin{equation}
	\ubex_i(A) := \inf\{q \in \R : \uBC_i^{\mathrm{var},q}(A) = 0\},
\end{equation}
where $\uBC_i^{\mathrm{var},q}(A)$
representing the critical scaling dimension for the $i$-th support measure encoded in the leading asymptotic term of the corresponding basic function. 
These exponents now provide a refined dimensional hierarchy that captures both the overall fractal dimension and the contributions of different geometric features to this dimension which is stated in the following theorem.

\begin{theorem}[Properties of basic scaling exponents \cite{RaWi1}] \label{thm:basic-exp}
    Let $A\subseteq\R^d$ be a nonempty compact set.
 For each $i\in I_d$ one of the following is true:
    \begin{enumerate}
        \item[{\rm (a)}]  $\mu_i(A;\cdot)\equiv 0$ and so, by convention, we let $\un{\frk{m}}_i(A)=\ov{\frk{m}}_i(A):=-\infty$;
        \item[{\rm (b)}] $i\leq \lbex_i(A)\leq \ubex_i(A)$.
        %The following inequalities hold:
      %  \begin{equation} \label{eq:mi-ineq}
	%	\leq \udimout_M A \quad  \text{ and }
	%\end{equation}	
    \end{enumerate}
       It is always true that  $\mu_0(A;\cdot)\not\equiv 0$. Furthermore,
     \begin{equation}\label{eq:dim_max_m}
     \max_{i\in I_d} \ubex_i(A) =\udimout_M A \quad \text{ and } \quad \max_{i\in I_d}\lbex_i(A)\leq \ldimout_M(A).
 \end{equation}
\end{theorem}

Recall that the notion of upper {\em surface} or {\em $S$-dimension} was introduced in \cite{RaWi2010} by $\udim_S A:=\inf\{s:\uS^s(A)=0\}$, where
    %\edz{Steffen, in our paper we use the normalized $d-1$ dimensional Hausdorff measure so please double check if it is the same in your paper with Jan, otherwise maybe the constant in the formula for the S-content needs to be adjusted or removed? No problem. The same normalization is used in RaWi2010}.
     \begin{align}\label{eq:S-content}
     \uS^s(A):=\limsup_{r\to 0^+} \frac{\sH^{d-1}(\partial A_r)}{\omega_{d-s} r^{d-1-s}}
     \end{align}
     is the $s$-dimensional \emph{upper $S$-content} of a compact set $A\subseteq\R^d$ and analogously for the lower conterpart.
     Note that by \cite[Corollary 3.6]{RaWi2010} one has $\udimout_M A=\udim_S A$ but, in general, one only has $\ldim_S A\leq\ldimout_MA$.

This result reveals a connection between the measure-theoretic structure encoded in support measures and the dimensional properties that characterize fractal sets. 

A direct implication of Theorem~\ref{thm:basic-exp} is that, for all indices $i$ satisfying $i>\ldimout_M A$, statement (a) is true (which means that $\mu_i(A;\cdot)\equiv 0$ and consequently $\lbex_i(A)=\ubex_i(A)=-\infty$).
Actually, a slightly stronger claim holds.
By \cite[Proposition 2.4]{Last06}, if $\sH^{k}(\partial A)=0$ for some nonempty closed set $A$, then $\mu_k(A;\cdot)\equiv 0$.
Since the Hausdorff dimension of $\partial A$ may be strictly less than the (lower) outer Minkowski dimension of $A$, $\mu_i(A;\cdot)$ may vanish for some more indices $i$ below $\ldimout_M A$, specifically for all $i>\lfloor{\dim}_H A\rfloor$. %. In particular, for any compact set $A\subseteq\R^d$,
%all the terms in the Steiner formula \eqref{eq:v-par} with an index $i>\lfloor{\dim}_H A\rfloor$ vanish.
However, a general upper bound for the basic scaling exponents is provided only by the (upper/lower) outer Minkowski dimension and not by the Hausdorff dimension, as equation \eqref{eq:dim_max_m} makes clear.

Indeed, for all indices $i\in I_d$ with $i\ne 0$, there exist sets $A$ for which $\mu_i(A;\cdot)\equiv 0$; that is, assertion (a) can occur for each such $i$. Additionally, the lower bound in part (b) is sharp for all $i\in I_d$, meaning equality can be attained. Compact convex sets $K\subseteq\R^d$ provide examples of both cases: if the affine dimension of $K$ is at least $i$, then $\bex_i(K)=i$; otherwise, $\bex_i(K)=-\infty$.

Moreover, if $K$ has affine dimension $k\le d-1$, then $\dimout_M K = \dim_M K = k = \max\{{\bex_i(K): i\in I_d}\}$. Such sets illustrate that $\ldimout_M K$ (and similarly $\udimout_M K$) gives a \emph{sharp} upper bound for all $\lbex_i(K)$, as stated in \eqref{eq:dim_max_m}. 
For the lower basic scaling exponents, only the upper bound in \eqref{eq:dim_max_m} is asserted here;
a corresponding equality need not hold in general.

%While \eqref{eq:dim_max_m} provides an exact value for the maximum upper scaling exponent, we conjecture that a corresponding equality does not hold in general for the lower scaling exponents. Typically, $\ldimout_M A$ only bounds $\max{\lbex_i(A)}$ from above, and $\ldim_S A$ from below, without equality.

\section{Examples}

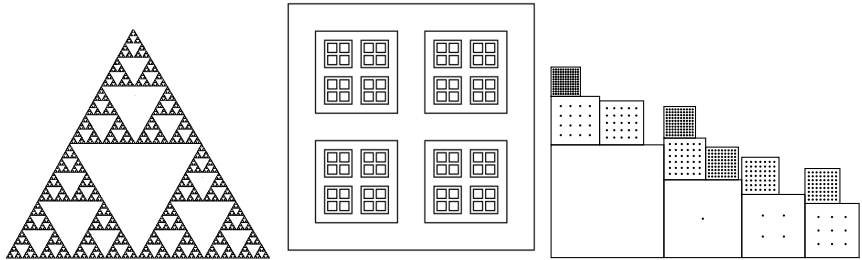
\begin{figure}[ht]
    \centering
    \includegraphics[width=0.3\textwidth]{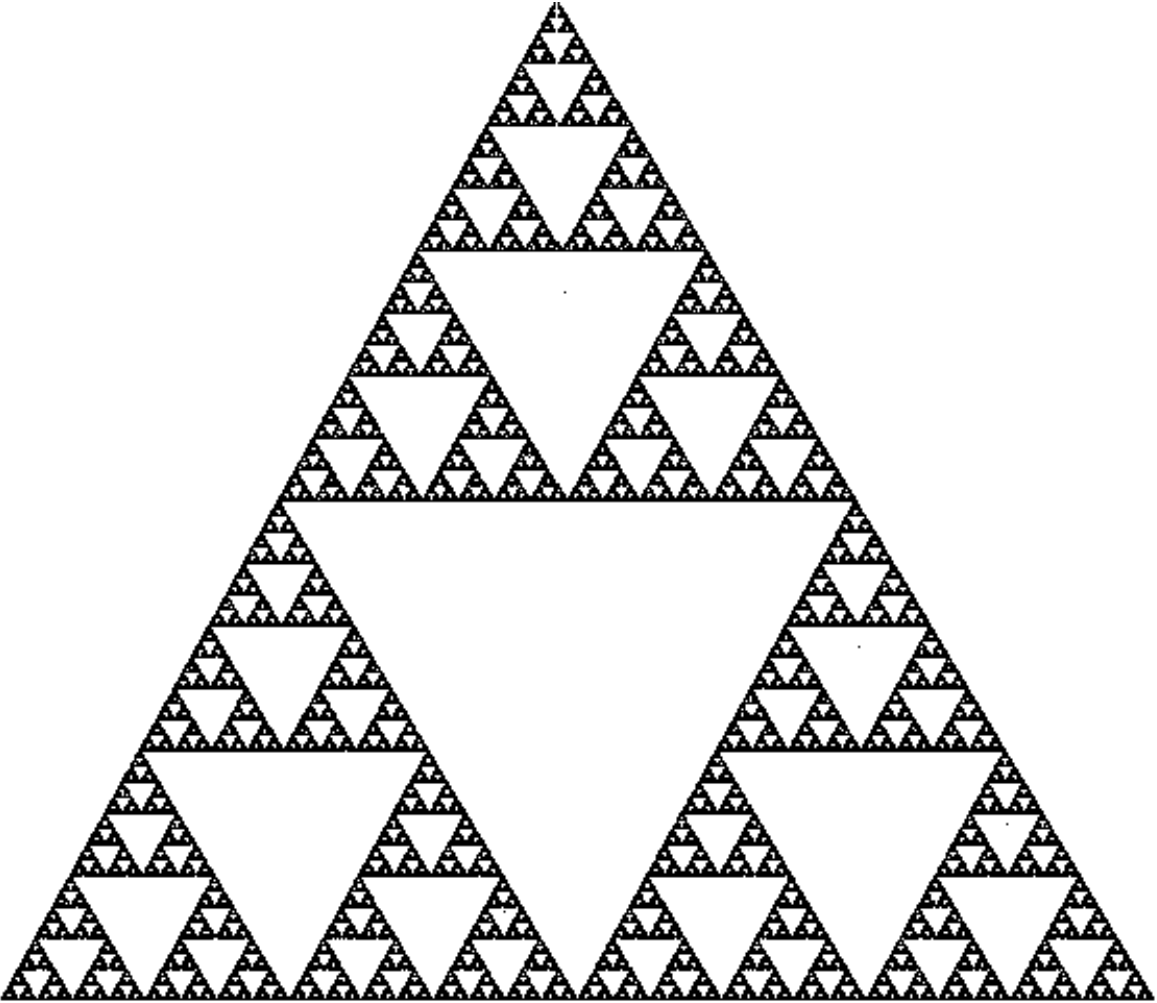}\ \ 
    \includegraphics[width=0.3\textwidth]{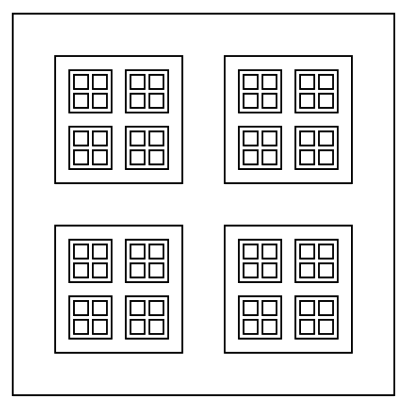}\ 
    \begin{tikzpicture}[scale=1.5]
	\draw (0, 0) rectangle (1,1);
		\equidistantPoints{xshift=1cm}{0.69}{1}{3};
		\equidistantPoints{xshift=1cm+0.69cm}{0.56}{2}{3}
		\equidistantPoints{xshift=1cm+0.69cm+0.56cm}{0.48}{3}{3}
		\equidistantPoints{yshift=1cm}{0.43}{4}{3}
		\equidistantPoints{yshift=1cm,xshift=0.43cm}{0.39}{5}{3}
		\equidistantPoints{yshift=0.69cm,xshift=1cm}{0.37}{6}{3}
		\equidistantPoints{yshift=0.56cm,xshift=1cm+0.69cm}{0.33}{7}{3}
		\equidistantPoints{yshift=0.48cm,xshift=1cm+0.69cm+0.56cm}{0.31}{8}{3}
		\equidistantPoints{yshift=0.69cm,xshift=1cm+0.37cm}{0.29}{9}{3}
		\equidistantPoints{yshift=0.69cm+0.37cm,xshift=1cm}{0.28}{10}{3}
		\equidistantPoints{yshift=1cm+0.43cm}{0.26}{11}{3}
	\end{tikzpicture}
   \caption{\label{fig:ex-joint} Examples of sets illustrating the possible relations between the basic exponents: $\bex_0<\bex_1$ (left, Sierpinski gasket), $\bex_0=\bex_1$ (middle, fractal window), $\bex_0>\bex_1$ (right, enclosed fractal dust), meaning that the support measure $\mu_1$ (left), $\mu_0$ (bottom) or both together (right) determine the outer Minkowski dimension.}
\end{figure}

Note that Theorem~\ref{thm:basic-exp} allows for any of the basic scaling exponents to dominate and thus determine the outer Minkowski dimension of a set. In $\R^2$, this concerns the two exponents $\bex_0$ and $\bex_1$. We now present three illustrative (classes of) examples demonstrating that all three cases: $\bex_0 < \bex_1$, $\bex_0 = \bex_1$, and $\bex_0 > \bex_1$ can occur; see also Figure~\ref{fig:ex-joint}.

More specifically, the Sierpinski gasket yields $\bex_0 = 0$ and $\bex_1 = \dimout_M > 1$, showing that the measure $\mu_1$ can encode the fractal behavior (Example~\ref{ex:SG}). In contrast, the \emph{fractal window} satisfies $\bex_0 = \bex_1 = \dimout_M > 1$ (Example~\ref{ex:FW}), while the \emph{enclosed fractal dust} exhibits $1 < \bex_1 < \bex_0 = \dimout_M$ (Example~\ref{ex:FD}). Thus, the outer Minkowski dimension can be determined by $\mu_1$, by $\mu_0$, or jointly by both.

In higher dimensions ($\R^d$, $d > 2$), similar behavior is expected: any of the exponents $\bex_i$ (or any combination thereof) may attain the maximum and thereby determine $\udimout_M$.

\begin{example}[The Sierpiński gasket] \label{ex:SG}
Let $SG$ denote the classical Sierpiński gasket constructed by iteratively removing $3^{n-1}$ open equilateral triangles of side length $2^{-n}$ from a unit equilateral triangle $A$, see Figure~\ref{fig:ex-joint}~(left). The limit of this process defines $SG$.
It is not too difficult to analyze the normal bundle $N(SG)$, which decomposes as $N(SG)=N(A)\cup N^{-}$, where $N^-$ consists of pairs $(x,u)$ with $x$ on the boundary of a removed triangle and $u$ the corresponding inward unit normal. For these points, the boundary is locally flat, implying $k_1(SG,x,u)=0$.

Note that $\mu_0(SG,\cdot)$ is supported on the vertices of $A$, and also equal to $\mu_0(A,\cdot)$ (see also \cite[Example 4.8]{HugLasWeil}), hence $\fb_0(SG;t)=1$ for all $t>0$. Thus, $\BC_0^0(SG)=1$ and $\bex_0(SG)=0$.

Since determining $\bex_1(SG)$ is much more involved, we only give a very rough sketch here for illustration purposes.
Note that the set $N(A)$ contributes a constant $\fb_1(A;t)=3/2$, representing half the boundary length of $A$. For $T \in \Delta_k$ (the set of triangles removed at step $k$), each has side length $2^{-k}$ and inradius $g\cdot 2^{-k+1}$, where $g=1/(4\sqrt{3})$. The contribution of $T$ to $\fb_1(SG;t)$ vanishes for $t \geq g \cdot 2^{-k+1}$, and for smaller $t$ it equals
\[
\ell(T,t) = 3(2^{-k-1} - \sqrt{3} t), \quad \text{for } t < g \cdot 2^{-k+1},
\]
see Figure \ref{fig:SG} for visualization.
\begin{figure}[ht]
    \centering
    \includegraphics[width=0.7\textwidth]{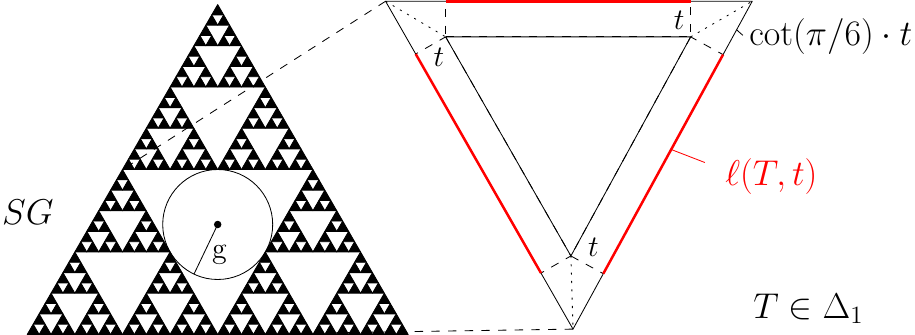}
    \caption{
    \label{fig:SG} The Sierpi\'nski gasket $SG$ from Example~\ref{ex:SG}. The largest removed triangle $T$ has inradius $g=1/(4\sqrt{3})$. The set of footpoints in $\pa T$ of pairs in the generalized normal bundle of $SG$ that contribute to $\beta_1(SG;t)$ is marked in red (and the contribution is given by $\ell(T,t)$).
    }
\end{figure}

Fix $n\in\N$ and suppose $t \in [g\cdot 2^{-n}, g\cdot 2^{-n+1})$. Then, all $\Delta_k$ for $1\leq k\leq n$ contribute, resulting in
\begin{equation} \label{eq:bf-SG}
\begin{aligned}
\fb_1(SG;t) &= \fb_1(A;t) + \sum_{k=1}^n \sum_{T\in\Delta_k} \ell(T,t) = \tfrac32+\sum_{k=1}^n 3^{k-1} \cdot 3(2^{-k-1}-\sqrt{3} t)\\
    &=\left(\tfrac32\right)^{n+1} - \tfrac 32 \sqrt{3} (3^n-1)\cdot t.
\end{aligned}
\end{equation}

From this expression, upper and lower bounds on $\fb_1(SG;t)$ are not difficult to obtain (one bounds $n$ in terms of $t$) and yield that the correct  exponent is, as expected, equal to the well-known (outer) Minkowski dimension of the Sierpinski carpet  
% Using $t \geq g\cdot 2^{-(n+1)}$, we obtain
% \[
% \fb_1(SG;t) \leq 3\left( \tfrac{3}{2} \right)^n \leq 3c^{-1} t^{1 - \log_2 3}, \quad c = (g/2)^{1 - \log_2 3},
% \]
% so that $\ubex_1(SG) \leq \log_2 3$.
% 
% For the lower bound, a similar calculation gives:
% \[
% \fb_1(SG;t) > \left( \tfrac{3}{2} \right)^{n+2} - \tfrac{3}{2} \quad \text{and} \quad t^{\log_2 3 - 1} \fb_1(SG;t) > \tfrac{5}{4} g^{\log_2 3 - 1}.
% \]
% Thus $t^{D-1}\fb_1(SG;t)$ stays uniformly bounded below for small $t$, implying $\lbex_1(SG) \geq D := \log_2 3$. Altogether, we conclude:
$
\bex_1(SG) = \log_2 3=\dimout_M(SG).
$
Note, however, that the basic content $\BC_1^{\log_2 3}(SG)$ does not exist as a limit, as it is the case with the Minkowski content of $SG$.
\end{example}

The next example demonstrates the situation when both support measures, $\mu_0$ and $\mu_1$ determine the outer Minkowski dimension of the set.

\begin{example}[The fractal window] \label{ex:FW}
We consider a family of compact sets $A=A(r) \subseteq \R^2$, parametrized by $r \in (0,1/2)$.
For $r\in[1/4,1/2)$ both basic scaling exponents coincide and attain a prescribed value in the interval $[1,2)$. That is, for each $r\in[1/4,1/2)$, the exponents satisfy
\[
\bex_0(A) = \bex_1(A) = \dimout_M A.
\]
In this case, the outer Minkowski dimension of the set is captured equally by both support measures $\mu_0$ and $\mu_1$. The set $A$ arises as an inhomogeneous self-similar set generated by four similarity transformations in the plane, each with contraction factor $r$; see Figure \ref{fig:ex-joint} (middle).

We begin with the boundary $B$ of the unit square—referred to as the \emph{initial frame}—and apply four contracting similarities $\Phi_1,\ldots,\Phi_4$, each scaling by $r$ and translating by vectors $\mathbf{x}_1=(p,p)$, $\mathbf{x}_2=(p,2p+r)$, $\mathbf{x}_3=(2p+r,p)$, and $\mathbf{x}_4=(2p+r,2p+r)$, where $p := (1-2r)/3$. This ensures that the images $\Phi_i(B)$ are uniformly spaced and disjoint from one another and from $B$ by distance $p$. Iterating this process infinitely yields a unique nonempty compact set $A$ satisfying the inhomogeneous fixed-point equation:
\[
A = \bigcup_{i=1}^4 \Phi_i(A) \cup B.
\]

The associated homogeneous attractor forms a Cantor dust with strong separation, whose Minkowski, Hausdorff, and similarity dimensions all equal the unique solution $s = \log_{1/r}4$ of the Moran equation $1 - 4r^s = 0$; see \cite{Hut81}.
However, as an inhomogeneous self-similar set (see \cite{MR3881123,MR3024316}), $A$ has outer Minkowski dimension
$
\dimout_M A = \max\left\{1, \log_{1/r}4\right\}.
$

Omitting details, we discuss here the basic exponents of the set $A$.
If $r \in (1/4, 1/2)$, then both basic exponents of $A$ coincide with the outer Minkowski dimension:
\[
\bex_0(A) = \bex_1(A) = \log_{1/r} 4 = \dimout_M(A)>1.
\]
Thus, the outer Minkowski dimension is jointly determined by $\mu_0$ and $\mu_1$.
Intuitively, similarly as in the case of the Sierpinski gasket, the measure $\mu_0$ is supported on the vertices of all the copies of the initial frame $B$ as opposed to the Sierpinski gasket where it was supported only on the three vertices of the initial triangle.
On the other hand, the measure $\mu_1$ is supported on the sides of all the copies of the initial frame $B$.
Since both, the vertices and the sides proliferate by the same similarity law, it is not unexpected that both of the basic exponents are equal.

However, for $r \in (0, 1/4)$, the similarity dimension stemming from the Moran equation is below 1 which equals to the topological dimension of $A$ and is captured by the measure $\mu_1$.
Hence, in this case we have $\bex_0(A) = \log_{1/r} 4 < 1 = \bex_1(A) = \dimout_M(A)$. In this range, $A$ is Minkowski measurable and forms a rectifiable curve with
\[
\BC_1^1(A) = \text{length}(A) = \frac{4}{1 - 4r} = \tfrac{1}{2} \Minkout{1}(A).
\]

At the critical value $r = 1/4$, a phase transition occurs: although $\bex_0 = \bex_1 = \dimout_M(A) = 1$, the basic function $\fb_1(A;t)$ exhibits logarithmic terms in $t$ (see \cite{RaWi2}), and $\BC_1^1(A) = \infty$, reflecting the infinite length of the associated curve. 
One may view $r = 1/4$ as marking the transition between regular (finite-length) and fractal geometry.
Notably, the basic function $\fb_0$ captures the self-similar structure of the set $A$ across the entire range $r \in (0, 1/2)$. Even when $A$ is rectifiable ($r < 1/4$), it retains a form of `subcritical fractality' encoded in the exponent $\bex_0$.
\end{example}

In the next example, a family of compact sets $A \subseteq \R^2$ is described for which the leading term in the Steiner formula is entirely determined by the support measure $\mu_0(A,\cdot)$. That is, the Minkowski dimension equals the $0$-th basic exponent, while the 1st basic exponent is strictly smaller:
$
\bex_1(A) < \bex_0(A) = \dim_M A.
$
This illustrates that the basic exponents are not necessarily ordered and may not satisfy $\bex_0 \leq \bex_1 \leq \ldots$ in general. Identifying necessary and sufficient conditions for such monotonicity remains an open problem.

\begin{example}[Enclosed fractal dust] \label{ex:FD}
The {\em enclosed fractal dust} is built from square boundaries (or `frames') $S_j := \partial Q_j$, where each $Q_j \subseteq \R^2$ is a closed square of side length $\ell_j > 0$, and $\text{int}(Q_j) \cap \text{int}(Q_k) = \emptyset$ for $j\ne k$, forming a packing. Inside each $Q_j$, $n_j^2$ points are placed uniformly on a square grid see Figure \ref{fig:ex-joint} (right).

The side lengths $(\ell_j)_{j\in\N}$ form a fractal string $\mathcal{L}$ in the sense of \cite{LF12}, and we assume that
$
\sum_{j=1}^\infty \ell_j^2 < \infty,
$
ensuring that the total area of the union $\bigcup_j Q_j$ is finite and can be contained in some compact set $K \subseteq \R^2$. %The grid spacing for placing the $n_j^2$ points within $Q_j$ is chosen as
%\begin{equation} \label{eq:r_j}
%r_j := \frac{\ell_j}{n_j + 1}.
%\end{equation}

While the construction can be generalized (e.g.\ to Minkowski measurable or non-measurable cases by varying the sequences $(\ell_j)$ and $(n_j)$, cf.\ \cite{LF12}), we focus here on a specific setup ensuring bounded area and packability.

Let $A_j$ denote the union of the square boundary $S_j := \partial Q_j$ and the $n_j^2$ equidistant points placed inside $Q_j$. The \emph{enclosed fractal dust} $A$ is then defined as
\[
A := \overline{\bigcup_{j=1}^{\infty} A_j}.
\]

In general, the arrangement of the $A_j$ inside $\R^2$ can affect the geometric properties of $A$, including its Minkowski dimension and scaling exponents. To ensure well-controlled behavior, we specify $(\ell_j)$ as a \emph{Riemann-type fractal string}:
\[
\ell_j := j^{-\alpha}, \quad \text{with } \alpha \in \left(\tfrac{1}{2}, \tfrac{2}{3}\right],
\]
which, as shown in \cite{MR4227795}, allows for a perfect packing into a square $K$ of area $\zeta(2\alpha)$, where $\zeta$ is the Riemann zeta function. Thus, all $Q_j$ (and hence $A_j$) are contained in a compact square $K \subseteq \R^2$.

For a fixed integer $m \in \N$, we define the number of grid points in $Q_j$ by
\[
n_j := \ell_j^{-m/\alpha} - 1 = j^m - 1,
\]
so that the grid spacing becomes
$
r_j := \frac{\ell_j}{n_j + 1}=j^{-(\alpha+m)}.
$

% Using fractal zeta function techniques (cf.\ \cite{RaWi2}), explicit formulas for the basic functions $\fb_0(A;t)$ and $\fb_1(A;t)$ valid for sufficiently small $t > 0$ are obtained:
% \begin{align}
% \fb_0(A;t) &= a_1 t^{-\frac{1 + 2m}{\alpha + m}} + a_2 t^{-\frac{1 + m}{\alpha + m}} + a_3 t^{-\frac{1}{\alpha + m}} + a_4, \label{eq:fb0_A} \\
% \fb_1(A;t) &= b_1 t^{1 - \frac{1 + m}{\alpha + m}} + b_2 t^{1 - \frac{1}{\alpha + m}} + b_3 + b_4 t, \label{eq:fb1_A}
% \end{align}
% where the constants $a_i$, $b_i$ depend explicitly on $\alpha$ and $m$ and are computable.

The basic exponents of the set $A$ behave as follows:
\[
\bex_0(A) = \frac{1 + 2m}{\alpha + m} > \bex_1(A) = \frac{1 + m}{\alpha + m} = \bex_0(A) - \frac{m}{\alpha + m} > 1,
\]
showing that $\mu_0$ governs the leading term in the Steiner formula for $A$. Moreover, it can be shown that the corresponding basic contents $\BC_0^{\bex_0}(A)$ and $\BC_1^{\bex_1}(A)$  exist in $(0,+\infty)$; see \cite{RaWi2} for their explicit expressions.
Intuitively, one can explain this phenomenon by noting that the proliferation of points (on which the measure $\mu_0$ is supported) is chosen to be much ``faster'' than the proliferation of squares (with $\mu_1$ being supported on the sides of these squares).
Hence, the points govern the leading term in the Steiner formula and the fact that both basic contents exist can be intuitively explained by the fact that there is no self-similarity in the construction of the squares, nor in the way the points are placed, and hence there is no underlying oscillation in the behavior of the basic functions.

Note that by varying $\alpha \in (1/2, 2/3]$ and $m \in \N$, one can obtain $\bex_0(A)$ in the full range $[9/5, 2)$ and $\bex_1(A)$ in $(1, 4/3)$, with $\bex_0(A) > \bex_1(A)$ holding in all cases. This demonstrates that $\mu_0$ alone may dominate the asymptotic behavior, and provides concrete examples where the basic exponents are not monotonically increasing.
\end{example}

\section{Complex Dimensions and Fractal Zeta Functions}\label{sec:compdim}

In this section we briefly recall basic facts about the theory of fractal zeta functions and complex dimensions as developed in \cite{FZF,lapidus2017distance,zbMATH07043493,zbMATH06828565,zbMATH06723342,lapidus2018fractal} which generalizes to $\R^d$, $d\geq 2$ the one-dimensional theory of complex dimensions for fractal strings \cite{LF12}. The aim is to demonstrate how the basic functions and corresponding basic scaling exponents provide crucial groundwork for linking the measure-theoretic framework of support measures with the analytical theory of fractal zeta functions and complex dimensions. Hence, complex analytical methods can be applied to extract geometric information from support measures.

Direct analysis of the basic functions $\fb_i$ can be made with considerable effort for various classes of sets--like self-similar tilings, the fractal window, and the enclosed fractal dust--it often provides limited explicit information about the corresponding basic contents since one usually extracts only the leading asymptotic terms.

To overcome this and extract even more information from the basic functions, we now outline an alternative and elegant approach using the theory of complex dimensions and fractal zeta functions developed in \cite{FZF}. Specifically, for any nonempty compact set $A \subseteq \R^d$, the (fractal) {\em Lapidus zeta function}, also known as the {\em distance zeta function}, is defined by
\begin{equation} \label{eq:distance-zeta}
\zeta_A(s) := \int_{A_\e \setminus A} \!\! \dist(z, A)^{s-d} \, \mathrm{d}z, \quad s \in \C,
\end{equation}
for fixed $\e > 0$. This slightly modified version, adapted to the relative domain $\Omega := A_\e \setminus A$, avoids assuming $\mathcal{L}^d(A) = 0$ and better aligns with our focus on outer Minkowski dimension. In the terminology of \cite[Ch.~4]{FZF}, we are working with the relative fractal drum (RFD) $(A, \Omega)$.

The integral in \eqref{eq:distance-zeta} converges absolutely and defines a holomorphic function on the half-plane $\{\re s > \udimout_M A\}$. Under mild assumptions, namely, if $D = \dimout_M A$ exists and $\lMinkout{D}(A) > 0$, then $D$ is a pole of $\zeta_A$ provided it admits a meromorphic continuation; see \cite[Thm.~2.1.11 and 4.1.7]{FZF}.

The set of \emph{complex dimensions} of $A$ (with respect to a domain $W \supseteq \{\re s > \udimout_M A\}$) is then defined as the set of poles of this meromorphic extension. Notably, the complex dimensions are independent of the choice of $\e$, since changing $\e$ alters $\zeta_A$ only by a adding an entire function. These dimensions generalize the classical Minkowski dimension and admit geometric interpretations, including asymptotic tube formulas involving residues at the poles; see \cite[Ch.~5]{FZF} known as {\em fractal tube formulas} which are another generalization of the classical Steiner formula.
Namely, omitting the details, for compact sets $A\subset\R^d$, satisfying certain assumptions, the following fractal tube formula holds:
\begin{equation} 
 \label{point_form}\nonumber
V_d(A_{\e}\setminus A)=\sum_{w\in\mathcal{P}(\zeta_A)} \frac{\eps^{d-w}}{d-w}\res\left({\zeta}_{A}(s),w\right),
\end{equation}
where $\mathcal{P}(\zeta_A)\subseteq\C$ is the set of simple complex dimensions of $A$.

The main results that bridges the theory of complex dimension with the basic functions and the general Steiner formula of Hug, Last and Weil is the following. 

\begin{theorem}[Basic zeta functions for compact sets \cite{RaWi2}]
	\label{thm:shell-decomp}
	Let $A$ be a compact subset of $\R^d$. %such that $\udimout_MA<d$ and fix $\e>0$.
	Then for all $s\in\C$ such that $\re s>\udimout_MA$, the following functional equation holds:
	\begin{equation}
		\label{eq:eqf-shell}
		\z_A(s)=\sum_{i=0}^{d-1}\o_{d-i}\breve{\z}_{A,i}(s),		
	\end{equation}
	where the {\em $i$-th basic zeta function of} $A$, $\breve{\z}_{A,i}$, for $i\in I_d$ is defined as
	\begin{equation}
		\label{eq:support-zeta}
		\breve{\z}_{A,i}(s)=\breve{\z}_{A,i}(s;\e):=\int_0^{\e}t^{s-i-1}\fb_i(t)\d t.
	\end{equation}
	
	Furthermore, the integral defining $\breve{\z}_{A,i}$ is absolutely convergent, and hence, holomorphic, in the open half-plane $\{\re s>\ov{\frk{m}}_i\}$.
\end{theorem}

From \eqref{eq:support-zeta}, under appropriate hypotheses concerning the growth of the basic zeta function, one can reconstruct each $\fb_i(t)$ via a Mellin inversion process, similarly as in \cite[Ch.\ 5]{FZF}; see \cite{RaWi2}.
This leads not only to the exponent $\bex_i(A)$, but also to the corresponding basic contents $\BC_i^{\bex_i}$ (if it exists), as well as higher-order or oscillatory correction terms. 

Importantly, two types of functional equations derived in \cite{RaWi2} allow one to compute $\breve{\zeta}_{A,i}$ in closed form without, a priori, directly knowing $\fb_i$—highlighting the computational power of the zeta-function approach.

\begin{theorem}[Functional equations for basic zeta functions of a compact set \cite{RaWi2}]\label{thm:basic-dist-funct}
	Let $A$ be a compact subset of $\R^d$ and fix $\e>0$. Let $i\in I_d$.
	Then for all $s\in\C$ such that $\re s>\ov{\frk{m}}_i$ the following functional equations hold:
	\begin{equation}
		\label{eq:reach-zeta}
		\breve{\z}_{A,i}(s)=\frac{1}{s-i}\int_{N(A)}(\min\{\dl(A,x,u),\e\})^{s-i}\mu_i(A;\d (x,u)),
	\end{equation}
	and
 	\begin{equation}
		\label{eq:basic-dist}
		\breve{\z}_{A,i}(s)=\int_{A_{\eps}\setminus A}\dist(z,A)^{s-i-1}\breve{K}_i(z)\,\d z,
	\end{equation}
	where both of the integrals above are absolutely convergent and hence, holomorphic in the open half-plane $\{\re s>\ov{\frk{m}}_i\}$ while
	\begin{equation}\label{eq:breve_K}
		\breve{K}_i(z)=\frac{1}{\o_{d-i}}\prod_{m=1}^{d-1}\frac{1}{1+\dist(z,A)k_m(A,\Pi_A(z))}\sum_{\substack{I\subseteq\{1,\ldots,d-1\} \\ |I|=d-1-i}}\prod_{l\in I}k_l(A,\Pi_A(z))
	\end{equation}
	and the sum extends over all subsets of $\{1,\ldots,d-1\}$ of cardinality $d-1-i$.
\end{theorem}

\medskip

Note that in expression \eqref{eq:reach-zeta} one does not need to know, a priori, the basic function $\fb_i$ in order to compute the integral since this is circumvented by using the local reach function $\dl(A,\cdot,\cdot)$, which is usually more practical to deal with.
Neverthless, one can later reconstruct $\fb_i$ by using Mellin inversion techniques and knowledge about the poles of meromorhic extension of the basic zeta function $\breve{\z}_{A,i}$; see Equation \eqref{point_form_w} just below.
Even more practical is the functional equation \eqref{eq:basic-dist} where the integration is done in $\R^d$ which is usually much more convenient than integration over the generalized normal bundle $N(A)$.

Omitting the details for the sake of brevity (see \cite{RaWi2}), we give here roughly the reconstruction formula for the basic functions via the associated basic zeta functions.
Let $A$ be a compact subset of $\R^d$ and let $i\in I_d$ such that its $i$-th basic zeta function $\breve{\zeta}_{A,i}(\cdot;\e)$ is meromorphic in all of $\C$ and satisfies appropriate growth conditions\footnote{See strong languidity conditions in \cite[Ch.\ 5]{FZF} or \cite{RaWi2}.}.
Then, there exists $\e_0>0$ such that the following exact formula for ${\fb}_i(A;t)$ holds for all $t\in(0,\e_0)$$:$
\begin{equation}\label{point_form_w}
\fb_i(A;t)=\sum_{w\in\mathcal{P}(\breve{\zeta}_{A,i}(\cdot;\e),\C)}\res\left(t^{i-s}\breve{\zeta}_{A,i}(s;\e),w\right),
\end{equation}
where $\mathcal{P}(\breve{\zeta}_{A,i},\C)$ is the set of all the poles of $\breve{\zeta}_{A,i}$.
% and $\kappa$ is the exponent occurring in the statement of hypotheses {\bf L1} and {\bf L2'}.
To demonstrate Theorem \ref{thm:basic-dist-funct} and application of formula \eqref{point_form_w} we go back to the Sierpinski gasket from Example \ref{ex:SG}.

\begin{example}[The Sierpinski gasket revisitied]
Let $SG$ be the Sierpinski gasket from Example \ref{ex:SG} and fix $\e > g$, where $g = 1/(4\sqrt{3})$ is the \emph{inradius} of $SG$.
Hence, the parallel set $SG_\e$ fully contains the unit triangle $A$. The set difference $SG_\e \setminus SG$ consists of all open triangles $O_n$ removed in the construction, along with the outer region $SG_\e \setminus A$, i.e., the $\e$-neighborhood of $SG$ outside $A$.

It is easy to see that
\[
k_1(SG, \Pi_{SG}(z)) = 
\begin{cases}
\infty, & \text{if } \pi_{SG}(z) \text{ is a vertex of } A,\\
0, & \text{otherwise}.
\end{cases}
\]
Let $D$ denote the union of the three circular sectors of radius $\e$ centered at the vertices of $A$.
Hence, one can now deduce by Equation~\eqref{eq:breve_K}, that the functions $\breve{K}_0$ and $\breve{K}_1$ are given by
\[
\breve{K}_0(z) = 
\begin{cases}
\frac{\dist(z,SG)^{-1}}{\omega_2}, & z \in D,\\
0, & z \in SG_\e \setminus (SG \cup D),
\end{cases}
\quad
\breve{K}_1(z) = 
\begin{cases}
0, & z \in D,\\
\frac{1}{\omega_1}, & z \in SG_\e \setminus (SG \cup D).
\end{cases}
\]
Using now Equation~\eqref{eq:basic-dist}, one can easily compute:
\[
\begin{aligned}
\breve{\zeta}_{SG,0}(s) &= \int_{SG_\e \setminus SG} \dist(z,SG)^{s-1} \, \breve{K}_0(z) \, \mathrm{d}z 
= \frac{1}{\omega_2} \int_D \dist(z,SG)^{s-2} \, \mathrm{d}z = \frac{2\pi \e^s}{\omega_2 s}.
\end{aligned}
\]
Thus, the meromorphic extension to all of $\C$ of $\breve{\zeta}_{SG,0}$  is given by right hand side and has a single simple pole at $s = 0$.

Next, for $\breve{\zeta}_{SG,1}(s)$ we again apply Equation~\eqref{eq:basic-dist}:
\[
\begin{aligned}
\breve{\zeta}_{SG,1}(s) &= \frac{3\e^{s-1}}{\omega_1(s-1)} + \frac{1}{\omega_1} \sum_{n=1}^\infty 3^{n-1} \int_{O_n} \dist(z,SG)^{s-2} \, \mathrm{d}z,
\end{aligned}
\]
where the first expression corresponds to three rectangles of size $1 \times \e$, while
the second integral is a sum over pairwise disjoint open triangles $O_n$ and, as shown easily (see also \cite[Example 4.2.24]{FZF}), one obtains the follwing closed expression
\[
\breve{\zeta}_{SG,1}(s) = \frac{3 \e^{s-1}}{\omega_1(s-1)} + \frac{6 (\sqrt{3})^{1-s} \cdot 2^{-s}}{\omega_1 s(s-1)(2^s - 3)}.
\]
Thus the meromorphic extension to all of $\C$ of $\breve{\zeta}_{SG,1}$ has simple poles at $s = 0$ and at the set of
\[
s = \log_2 3 + \tfrac{2\pi \I k}{\log 2}, \quad k \in \Z,
\]
while the apparent pole at $s = 1$ is removable since the residue equals zero.

Theorem~\ref{thm:shell-decomp} now reconstructs the distance zeta function of the Sierpinski gasket, $\zeta_{SG}$, computed in \cite[Prop.~3.2.3]{FZF} which was then used in \cite[Ex.~5.5.12]{FZF} to derive the fractal tube formula for $SG$, expressing the volume $V(SG_t)$ in terms of the complex dimensions, i.e., the poles of $\zeta_{SG}$.
	
\end{example}

Applying formula \eqref{point_form_w} to the Sierpiński gasket $SG$, we obtain the explicit expression:
$$
\fb_1(SG;t) = t^{1 - \log_2 3} \cdot \frac{3\sqrt{3}}{\log 2} \sum_{k \in \Z} \frac{(4\sqrt{3})^{-\nu_k} \, \mathrm{e}^{-2\pi \I k \log_{2} t}}{\nu_k(\nu_k - 1)} + \frac{3\sqrt{3}}{2} t
$$
for $t \in (0, g)$, where $\nu_k := \log_2 3 + \tfrac{2\pi \I k}{\log 2}$. The multiplicatively periodic Fourier series converges absolutely and defines a continuous function on $(0, g)$.
As a result, we recover the basic exponents and basic contents of $SG$:
$$
\bex_0(SG)=0, \quad \bex_1(SG) = \log_2 3, \quad \mathcal{M}_0^0(SG)=1
$$
while  we note that $\BC_1^{\log_2 3}(SG)$ does not exist. 
However, the upper and lower basic contents $\uBC_1^{\log_2 3}(SG)$ and $\lBC_1^{\log_2 3}(SG)$ exist, and equal the maximum and minimum of the periodic function defined by the Fourier series.

Finally, substituting into the general Steiner formula \eqref{eq:v-par} yields the well-known fractal tube formula for the Sierpiński gasket, consistent with \cite[Example 5.4.12]{FZF}—expressing $V(SG_t)$ explicitly in terms of the complex dimensions (i.e., the poles of $\zeta_{SG}$).

\section{Support Contents and the Parallel Set Perspective}

While basic functions provide direct access to the geometric information encoded in support measures, the theory also develops a parallel approach through support contents, which are based on the behavior of support measures for parallel sets. This dual perspective proves essential for connecting the measure-theoretic framework to classical approaches in fractal geometry.

For any compact set $A \subseteq \R^d$ and any index $i \in I_d$, the {\em $i$-th tube-support function} is defined as the mapping $\e \mapsto \mu_i(A_\e)$, which describes how the total mass of the $i$-th support measure changes as the set is dilated by parallel expansion. This function captures the global behavior of support measures under scaling, complementing the local information provided by basic functions.

The {\em $q$-dimensional upper $i$-th support content} is defined as
$$\uSC_i^{q}(A) := \limsup_{\e \to 0^+} \e^{q-i} \mu_i(A_\e),$$
with the upper $i$-th support scaling exponent given by
$$\usex_i(A) := \inf\{q \in \R : \uSC_i^{\mathrm{var},q}(A) = 0\},$$
where one uses the total variation analog of $\uSC_i^{q}(A)$ and the lower analogs are then defined in a standard way via the lower limit.  

The relationship between basic functions and support measures for parallel sets is established through a fundamental Steiner-type formula from \cite{HugLasWeil}. Namely, for any $\e > 0$ and any index $i \in I_d$, the support measures of parallel sets satisfy
$$\mu_i(A_\e) = \sum_{j=0}^i c_{i,j} \e^{i-j} \fb_j(\e),$$
where $c_{i,j}$ are universal positive constants that depend only on the dimension and indices involved.

This formula reveals the deep connection between the local geometric information encoded in basic functions and the global scaling behavior captured by support measures of parallel sets. The triangular structure of the coefficient matrix ensures that the relationship can be inverted, allowing for the expression of basic functions in terms of support measures of parallel sets.

The support scaling exponents exhibit important relationships with their basic counterparts summarized in the next theorem.

\begin{theorem}[Properties of support scaling exponents \cite{RaWi1}] \label{thm:support-exp}
    Let $A\subseteq\R^d$ be a nonempty compact set.
 For each $i\in I_d$ and each $\e>0$,  one has $\mu_i(A_\eps;\cdot)\not\equiv 0$,
        and
        \begin{equation}\label{eq:s_i-ineq}
		0\leq \lsex_i(A)\leq \usex_i(A)\leq \usex_{d-1}(A) =\udimout_M A.
  %\edz{Can we have $i\leq\lsex_i$ if $i\leq\lfloor\dim_HA\rfloor$?}
	\end{equation}
    Furthermore,
     \begin{equation}\label{eq:s<max}
     \usex_i(A)= \max_{j\in\{0,\ldots, i\}}\ubex_j(A).
 \end{equation}
 Similarly,
  \begin{equation}\label{eq:m<smax}
     \ubex_i(A)\leq \max_{j\in\{0,\ldots, i\}}\usex_j(A).
 \end{equation}
 Moreover,
 \begin{equation} \label{eq:sd-1-dim}
     \max_{i\in I_d}\usex_i(A)=\udimout_M(A)
     \quad \text{ and } \quad
      \ldim_S(A)\leq\max_{i\in I_d}\lsex_i(A)\leq \ldimout_M(A).
 \end{equation}
\end{theorem}

These relations show how the upper support scaling exponents aggregate the informa-
tion recorded separately by the upper basic scaling exponents.

It can be shown that in $\R^2$, the exponents $\usex_0$ and $\usex_1$ are always uniquely characterized by $\ubex_0$ and $\ubex_1$. Namely, one always has that $\usex_0(A)=\ubex_0(A)$, while 
$\usex_1(A)=\udimout_M A=\max\{\ubex_0(A),\ubex_1(A)\}$. Furthermore, the inequality $\usex_0(A)\leq\usex_1(A)$ is satisfied for every compact set $A\subseteq\R^2$. This ensures that $\usex_0$ cannot exceed $\usex_1$, which stands in contrast to what occurs with the basic exponents.

The inequality in \eqref{eq:m<smax} may indeed be strict, as is evident. Moreover, \eqref{eq:s<max}  shows that in every dimension
the upper support scaling exponents are nondecreasing, satisfying
\begin{align}
    \label{eq:increasing-exp}
    \usex_0(A)\leq\usex_1(A)\leq\cdots\leq\usex_{d-2}(A)\leq\usex_{d-1}(A).
\end{align}
Nevertheless, for spaces $\R^d$ where $d\geq 3$, whether \eqref{eq:increasing-exp} holds universally continues to be an unresolved open problem. %\edz{What about lower exponents!}

\section{The Bridge to Minkowski Content}

We comment here on the connection to classical concepts in fractal geometry, particularly the Minkowski dimension and content. This connection is most clearly manifested in the case of the $(d-1)$-th support measure, which corresponds to a generalized surface area measure.

For any nonempty compact set $A \subseteq \R^d$, the $(d-1)$-th support scaling exponent satisfies $\usex_{d-1}(A) = \udimout_M A$ by Theorem \ref{thm:support-exp}. Moreover, when the support content $\SC_{d-1}^s(A)$ exists for some $s \in [0,d)$, the outer Minkowski content exists and satisfies the explicit relationship
$$\Minkout{s}(A) = \frac{2}{d-s} \SC_{d-1}^s(A).$$

This connection extends far beyond the special case of surface measures. When basic contents $\BC_j^q(A)$ exist in $[0,\infty)$ for some $q\in[0,d)$ and for all indices $j \in I_d$, the outer Minkowski content admits a complete decomposition:
$$\Minkout{q}(A) = \frac{2}{d-q} \sum_{j=0}^{d-1} c_{d-1,j} \BC_j^q(A).$$

This decomposition reveals how different geometric features contribute to the overall Minkowski content, with each summand representing the contribution of "$j$-dimensional features" to the fractal content provideing insight into the geometric structure of fractal sets by showing how the overall dimensional properties emerge from the interplay of features at different dimensional scales.

\begin{acknowledgement}
The research of Goran Radunovi\'c was supported by the Croatian Science Foundation (HRZZ) grant IP-2022-10-9820 and by the Horizon grant 101183111-DSYREKI-HORIZON-MSCA-2023-SE-01.
\end{acknowledgement}
\ethics{Competing Interests}{The author has no conflicts of interest to declare that are relevant to the content of this chapter.}

\bibliographystyle{abbrv}
\bibliography{GR-bib}

\end{document}